\documentclass[twoside, a4paper,11pt]{article}
\usepackage{amsmath}
\usepackage{amssymb}
\usepackage{latexsym,a4wide}
\usepackage{xy}
\usepackage{hyperref}

\input xypic
\input xy
\xyoption{all}

\newtheorem{example}{Example}[section]
\newtheorem{definition}[example]{Definition}
\newtheorem{prop}[example]{Proposition}
\newtheorem{thm}[example]{Theorem}

\newtheorem{ex}[example]{Example}

\newenvironment{pf}{\noindent \textbf{Proof:}}{\rule{0em}{1ex}\hfill$\Box$\mbox{}}
\def\Ob{\operatorname{Ob}}

\begin{document}

\date{}

\title{QUADRATIC MODULES BIFIBRED OVER NIL(2)-MODULES}
\author{Hasan Atik and Erdal Ulualan}
\maketitle

\begin{abstract}
In this work, we show that the forgetful functor from the category of Baues's quadratic modules to that of nil(2)-modules is a bifibration.
\end{abstract}

\section*{Introduction}
Quadratic modules were introduced by Baues, \cite{baus}, as algebraic models for homotopy connected 3-types. They are related to simplicial groups,
crossed squares \cite{walery, loday} and  2-crossed modules \cite{con} as explained in \cite{AU1, baus}.
The freeness conditions for quadratic modules of groups in terms of free simplicial groups were also described in \cite{au}.

Brown and Sivera \cite{brown6} have shown that the forgetful functor $\Phi_1:\mathsf{XMod}\rightarrow  \mathsf{Gpd}$
 from the category of crossed modules of \textit{groupoids} to  the
category of\textit{ groupoids}, which sends a crossed module $M\rightarrow P $ to its base groupoid $P$, is a fibration and also a cofibration of categories.
 This result in the case of crossed modules of groups was appeared in \cite{brown4} and pursued in \cite{b-w1, b-w2}. This allows detailed computations of non-Abelian information on second relative groups. Analogous constructions in the crossed modules category in Lie algebras and
commutative algebras are given in \cite{casas} and \cite{porter1}, respectively.
For further accounts  of  fibred and cofibred categories and introduction to their literature, see \cite{brown9, gr, ts} and the references there.

In this paper, to shed some light on the quadratic module structure, we extend the result of Brown and
Sivera to the forgetful functor $\Phi:\mathsf{Quad}\rightarrow \mathsf{Nil}(2)$ from the category of quadratic modules to that of nil(2)-modules.
Thus, the main aim of this paper is  to show that quadratic modules are fibred and cofibred over
nil(2)-modules. This gives the notions of  \textit{pullback} and \textit{induced} quadratic module by a morphism of nil(2)-modules.
Clearly, this result and general methods given for quadratic modules are 2-dimensional analogues of the description of the
 \textit{pullback} and \textit{induced}  crossed modules  given in \cite{brown4, brown5,  brown9, brown6, casas, porter1}.

\section{Basic Definitions}

\subsection{Fibration and Cofibration of Categories}

The definitions and properties related to fibrations and cofibrations of
categories, some of which are classical, are recalled from \cite{brown6}.

\begin{definition}
Let $\Phi :\mathsf{X}\rightarrow \mathsf{B}$ be a functor. A morphism $\varphi :Y\rightarrow X$ in $\mathsf{X}$ over $u:=\Phi (\varphi )$ is called
cartesian if and only if for all $\upsilon :K\rightarrow J$ in $\mathsf{B}$
and $\theta :Z\rightarrow X$ with $\Phi (\theta )=u\upsilon $ there is a
unique morphism $\psi :Z\rightarrow Y$ with $\Phi (\psi )=\upsilon $ and $\theta =\varphi \psi .$

This is illustrated as follows.
\begin{equation*}
\xymatrix{Z\ar@{.>}[r]_\psi\ar@/^1.5pc/[rr]^\theta
&Y\ar[r]_\varphi&X&\ar[d]^{\mathsf{\Phi}}\\
K\ar@/^1.5pc/[rr]^{u\upsilon}\ar[r]_{\upsilon}&J\ar[r]_u&I&}
\end{equation*}
\end{definition}

A morphism $\alpha:Z\rightarrow Y$ is called vertical (with respect to $\Phi$) if and only if $\Phi(\alpha)$ is an identity isomorphism in $\mathsf{B}$.
In particular, for $I\in \mathsf{B}$ we write $\mathsf{X}_I$, called the
fibre over $I$, for the subcategory of $\mathsf{X}$ consisting of those
morphisms $\alpha$ with $\Phi(\alpha)=id_I,$

\begin{definition}
The functor $\Phi:\mathsf{X}\rightarrow \mathsf{B}$ is fibration or category
fibred over $\mathsf{B}$ if and only if for all $u:J\rightarrow I$ in $\mathsf{B}$ and $X\in \mathsf{X}_I$ there is a cartesian morphism $\varphi:Y%
\rightarrow X$ over u: such a $\varphi$ is called a cartesian lifting of X
along u.
\end{definition}

If $\Phi :\mathsf{X}\rightarrow \mathsf{B}$ is a fibration, then using the
axiom of choice for classes we may select for every $u:J\rightarrow I$ in
$\mathsf{B}$ and $X\in \mathsf{X}_{I}$ a cartesian lifting of X along u
\begin{equation*}
u_{X}:u^{\ast }X\rightarrow X.
\end{equation*}

Such a choice of cartesian lifting is called a cleavage or splitting of $\Phi$. If we fix the morphism $u:J\rightarrow I$ in $\mathsf{B}$, the
splitting gives a so-called reindexing functor
\begin{equation*}
u^{\ast }:\mathsf{X}_{I}\rightarrow \mathsf{X}_{J}
\end{equation*}
defined on objects by $X\mapsto u^{\ast }X.$ We can use this re-indexing
functor to get an adjoint situation for each $u:J\rightarrow I$ in $\mathsf{B}$.
The functor $u^*$ is said to give the objects \textit{pullback} by $u$.

\begin{definition}
Let $\Phi :\mathsf{X}\rightarrow \mathsf{B}$ be a functor. A morphism $\psi
:Z\rightarrow Y$ in $\mathsf{X}$ over $\upsilon :=\Phi (\psi )$ is called
cocartesian if and only if for all $u:J\rightarrow I$ in $\mathsf{B}$ and $\theta :Z\rightarrow X$ with $\Phi (\theta )=u\upsilon $ there is a unique
morphism $\varphi :Y\rightarrow X$ with $\Phi (\varphi )=u$ and $\theta
=\varphi \psi .$

This is illustrated as follows.
\begin{equation*}
\xymatrix{Z\ar[r]_\psi \ar@/^1.5pc/[rr]^\theta
&Y\ar@{.>}[r]_\varphi&X&\ar[d]^{\mathsf{\Phi}}\\
K\ar@/^1.5pc/[rr]^{u\upsilon}\ar[r]_{\upsilon}&J\ar[r]_u&I&}
\end{equation*}
\end{definition}

\begin{prop}\label{2}(\emph{\cite{brown6}}, Proposition 3.7)
Let $\Phi :\mathsf{X}\rightarrow \mathsf{B}$ be a fibration of categories.
Then $\psi :Z\rightarrow Y$ in $\mathsf{X}$ over $\upsilon :K\rightarrow J$
in $\mathsf{B}$ is cocartesian if and only if for all $\theta ^{\prime
}:Z\rightarrow X^{\prime }$ over $\upsilon $ there is a unique morphism $\psi ^{\prime }:Y\rightarrow X^{\prime }$ in $\mathsf{X}_{J}$ with $\theta ^{\prime
}=\psi ^{\prime }\psi $.\bigskip
\end{prop}
The functor $\Phi :\mathsf{X}\rightarrow \mathsf{B}$ is cofibration or
category cofibred over $\mathsf{B}$ if and only if for all $\upsilon
:K\rightarrow J$ in $\mathsf{B}$ and $Z\in \mathsf{X}_{K}$ there is a
cocartesian morphism $\psi :Z\rightarrow Z^{\prime }$ over $\upsilon :$ such a
$\psi $ is called a cocartesian lifting of $X$ along $\upsilon $.

If $\Phi:\mathsf{X} \rightarrow \mathsf{B}$ is a cofibration, for every morphism $v:K\rightarrow  J$ in $\mathsf{B}$ and an object $Z\in \mathsf{X}_K$, a cocartesian lifting of $Z$
$$
v_Z:Z\rightarrow  v_*(Z)
$$
along $v$ can be selected. Under these conditions, the functor $v_*$ is said to give the objects \textit{induced} by $v$.

\begin{prop}
Let $\Phi :\mathsf{X}\rightarrow \mathsf{B}$ be a functor that has a left
adjoint $D$. Then for each $K\in \Ob\mathsf{B},D(K)$ is initial in $\mathsf{X}_{K}$. In
fact if $\upsilon :K\rightarrow J$ in $\mathsf{B}$, then for any $X\in \mathsf{X}_{J}$ there is a unique morphism $\epsilon _{K}:DK\rightarrow X$ over $\upsilon $.
\end{prop}
\begin{thm}\label{3}(\emph{\cite{brown6}}, Theorem 4.2)
Let $\Phi :\mathsf{X}\rightarrow \mathsf{B}$ be a fibration of categories
which has a left adjoint D. Suppose that $\mathsf{X}$ admits pushouts. Let $\upsilon
:K\rightarrow J$ be a morphism in $\mathsf{B}$, and let $Z\in \mathsf{X}_{K}.$ Then a
cocartesian lifting $\psi :Z\rightarrow Y$ of $\upsilon $ is given precisely
by the pushout in $\mathsf{X}$:
\begin{equation*}
\xymatrix{D(K)\ar[rr]^{D(\upsilon)}\ar[dd]_{\epsilon_K}&&D(J)\ar[dd]^{\epsilon_J}\\ &&&\\ Z\ar[rr]_\psi&&Y}
\end{equation*}
\end{thm}

\section{Nil(2)-Modules Fibred and Cofibred Over Groups}
Throughout this paper all actions will be right. The left actions in some references
 will be rewritten by using right actions.

Recall from \cite{porter2} that a \textit{ pre-crossed module} is a group homomorphism $\partial
:M\rightarrow Q$ together with an action of $Q$ on $M$, written $m^{q}$ for $q\in Q$ and $m\in M$, satisfying the condition $\partial
(m^{q})=q^{-1}\partial (m)q$ for all $m\in M$ and $q\in Q$.

\textit{A nil(2)-module} (cf. \cite{baus}) is a pre-crossed module $\partial
:M\rightarrow Q$ with an additional `nilpotency' condition. This condition
is $P_{3}(\partial )=1$, where $P_{3}(\partial )$ is the subgroup of $M$
generated by Peiffer commutators $\left\langle
x_{1},x_{2},x_{3}\right\rangle $ of length $3$ for $x_i\in M$. The Peiffer
commutator in a pre-crossed module $\partial :M\rightarrow Q$ is defined by
\begin{equation*}
\left\langle x,y\right\rangle =x^{-1}y^{-1}x(y)^{\partial x}
\end{equation*}
for $x,y\in M$. For a pre-crossed module $\partial:M\rightarrow Q$, if $\langle M, M \rangle=1$, then it is called a \emph{crossed module} introduced by Whitehead (cf. \cite{wayted}).

A morphism $(g,f):(M\rightarrow Q) \longrightarrow (M'\rightarrow Q' )$
between  nil(2)-modules  is a pair of homomorphisms of groups $g:M\rightarrow M^{\prime }$ and $f:Q\rightarrow
Q^{\prime }$ such that $f\partial =\partial ^{\prime }g$ and the actions are
preserved, i.e. $g(m^{q})=g(m)^{f(q)}$ for any $m\in M,q\in Q$. This defines the category $\mathsf{Nil}(2)$ having nil(2)-modules as objects.

We have a forgetful functor from the category of nil(2)-modules to the category of groups
$$
\Phi_1:\mathsf{Nil}(2)\rightarrow \mathsf{Grp}
$$
which sends a nil(2)-module $(M\rightarrow Q)$ to the group $Q$.

This functor has also a left adjoint functor
$$
D_1:\mathsf{Grp}\rightarrow \mathsf{Nil}(2)
$$
assigns to a group $Q$ the trivial nil(2)-module $\{1\}\rightarrow Q$.

The following result in the cases of crossed modules of groups and groupoids appeared in \cite{brown4} and \cite{brown6} respectively,
described in terms of the crossed module of group(oid)s $\sigma_*M\rightarrow Q$ induced
from the crossed module of group(oid)s $M\rightarrow P$ by a morphism $\sigma:P\rightarrow Q$ in the category of group(oid)s. A similar result for nil(2)-modules over groups
 can be given as follows:
\begin{prop}
The forgetful functor $\Phi_1:\mathsf{Nil}(2)\rightarrow \mathsf{Grp}$ is fibred and cofibred.
\end{prop}
\begin{pf}
We give the pullback construction of nil(2)-modules to prove that $\Phi_1$ is fibred. Suppose that $\partial :M\rightarrow Q$ is a nil(2)-module and $\sigma
:P\rightarrow Q$ is a homomorphism of groups. Define $\sigma^*(M)$ to be the subgroup of $P\times M$ of elements $(p,m)$ such that $\sigma(p)=\partial(m).$
Let $\sigma _{1}:(p,m)\mapsto m$ and $\beta _{1}:(p,m)\mapsto p$. We obtain the following diagram
\begin{equation*}
\xymatrix{\sigma^*(M)\ar[d]_{\beta_{1}}\ar[r]^-{\sigma_{1}}&M\ar[d]^\partial
\\P\ar[r]_\sigma&Q.}
\end{equation*}
The action of $p^{\prime }\in P$ on $(p,m)\in
\sigma^*(M)$ can be given by $(p,m)^{p^{\prime }}=(p^{\prime -1}pp^{\prime},m^{\sigma(p')})$ and $\beta_1$ becomes  a nil(2)-module. Thus we have a nil(2)-module
with the base $P$.  This nil(2)-module $\beta_1$ is called the  \textit{pullback} of nil(2)-module $\partial$ along $\sigma$.
In the above diagram,  $(\sigma_1,\sigma)$ becomes  \textit{cartesian} morphism in the category of nil(2)-modules over the morphism  $\Phi_1(\sigma_1,\sigma)=\sigma$.
Therefore the forgetful functor $\Phi_1$ is a fibration of categories.

By direct construction, we prove that $\Phi_1$ is cofibred. To get it, we give the induced construction of nil(2)-modules.

  Let $\mu : M \rightarrow P$ be a nil(2)-module and let $f: P \rightarrow Q$ be a homomorphism of groups. We define $f_*(M) = F(M\times Q)$ to be the free group
on the set of elements $(m,q)$. The action of $q'$ on $(m,q)$ is $(m,q)^{q'}=(m,qq')$. We define as usual $\partial':F(M\times  Q)\rightarrow Q$ to be $(m,q)\mapsto
q^{-1}f(\mu m)q$. It is well known from \cite{brown9} that this gives a pre-crossed module over $F(M\times Q)\rightarrow Q$ with a map $i:M\rightarrow F(M\times Q)$ given
by $m\mapsto (m,1)$.

To make $i:M\rightarrow F(M\times Q)$ an operator morphism, we need factor $F(M\times Q)$ out by the relations $(m,q)(m',q)=(mm',q)$ and $(m^p,q)=(m,f(p)q)$ for $þp\in P$.
To make $\partial':F(M\times Q)\rightarrow Q$ a nil(2)-module involves factoring out triple Peiffer elements.

On generators, we have
$$\theta (m^{p})= (m^{p},1) = (m,f(p)1)
=((m,1))^{f(p)} = \theta (m)^{f(p)}$$ and $\partial'\theta (m)=\partial'((m,1))=f\mu (m)$, then we have the following diagram
\begin{equation*}
\xymatrix{&&N\ar@/^1.1pc/[ddl]^\upsilon\\
M\ar[d]_\mu\ar@/^1.1pc/[urr]^h\ar[r]^-\theta&F(M\times
Q)\ar@{.>}[ur]_{h'}\ar[d]^{\partial'}&\\ P\ar[r]_f&Q&}
\end{equation*}
in which $(h,f):(M\rightarrow P)\rightarrow (N\rightarrow Q)$ is a
nil(2)-module morphism and $h':{F(M\times Q)}\rightarrow N$ given by $(m,q)\mapsto h(m)^{q}$ is the necessary unique morphism in $\mathsf{Nil}(2)_{Q}$ which is the category
 of nil(2)-modules over the same group $Q$, and where $\theta$ is induced from the map $i$. Using Proposition \ref{2}, we get a nil(2)-module morphism
$(\theta,f):(M,P)\rightarrow (F(M\times Q),Q)$ which is cocartesian in $\mathsf{Nil}(2)$ over the
 homomorphism of groups $\Phi_1(\theta,f)=f$. Thus $\Phi_1$ is also a cofibration.
\end{pf}

\section{Quadratic Modules Fibred  over Nil(2)-Modules}
Quadratic modules were introduced by Baues in \cite{baus} as algebraic
models for topological $3$-types (i.e. homotopy connected 3-types).
\begin{definition}
A quadratic module $(\omega ,\delta, \partial)$ is defined as a diagram of $N $-groups
\begin{equation*}
\xymatrix{& C\otimes C\ar[d]^{w}\ar[dl]_{\omega} &\\ L \ar[r]_{\delta} &
M\ar[r]_{\partial}&N }
\end{equation*}
where $\partial$ is a nil(2)-module, $q:M\twoheadrightarrow(M^{cr})^{ab}=C$
is the natural projection map, $C\otimes C$ has the diagonal action, and the
following equalities hold

$\mathsf{QM1}.$ $\partial \delta=1$,

$\mathsf{QM2}.$ $\delta \omega (q(x)\otimes q(y))=w(q(x)\otimes q(y))=\langle x,
y\rangle$,

$\mathsf{QM3}.$ $\omega (q(\delta a)\otimes q(x))(q(x)\otimes q(\delta
a))=a^{-1}a^{\partial x},$

$\mathsf{QM4}.$ $\omega (q(\delta a)\otimes q(\delta b))=[a,b],$\newline
for $a,b\in L,x,y\in M$. The natural projection map $q:M\twoheadrightarrow(M^{cr})^{ab}=C$
is denoted on elements by $x\mapsto q(x)=\{x\}$ for $x\in M$.
\end{definition}

A map $\varphi :(\omega ,\delta ,\partial )\rightarrow (\omega ^{\prime
},\delta ^{\prime },\partial ^{\prime })$ between quadratic modules is given
by a commutative diagram, $\varphi =(l,m,n)$
\begin{equation*}
\xymatrix{ C\otimes C\ar[d]_{\varphi_{\ast}\otimes\varphi_{\ast}}
\ar[r]^-{\omega} & L \ar[d]_{l} \ar[r]^{\delta} & M \ar[d]_{m}
\ar[r]^{\partial} & N \ar[d]_{n} \\ C^{\prime }\otimes C^{\prime }
\ar[r]_-{\omega^{\prime}} & L^{\prime} \ar[r]_{\delta^{\prime}} & M^{\prime}
\ar[r]_{\partial^{\prime}} & N^{\prime} }
\end{equation*}
where $(m,n)$ is a morphism between pre-crossed modules which induces $\varphi _{\ast }:C\rightarrow C^{\prime }$ and where $l$ is an $n$-equivariant
 homomorphism. Let $\mathsf{Quad}$ be the category of quadratic
modules and of maps as in the above diagram.

There is a forgetful functor
 $$\Phi: \mathsf{Quad} \rightarrow \mathsf{Nil}(2)$$
 from the category of quadratic modules
to the category  of nil(2)-modules which sends a quadratic module
\begin{equation*}
\xymatrix{& C\otimes C\ar[d]^{w}\ar[dl]_{\omega} &&&\\ L \ar[r]_{\delta} &
M\ar[r]_{\partial}&N}
\end{equation*}
to its base nil(2)-module $\xymatrix{(M\ar[r]^{\partial}&N)}$.

The left adjoint of this functor is defined as follows. Recall from \cite{baus} that any nil(2)-module $\partial :M\rightarrow N$ yields a quadratic module
$\overline{\partial }:(1,w,\partial )$ given by
\begin{equation*}
\xymatrix{ & C\otimes C\ar[dl]_{1} \ar[d]^{w} \\ L \ar[r]_{w} & M
\ar[r]_{\partial} & N }
\end{equation*}
where $L=C\otimes C$. This quadratic module is called the quadratic module
associated to the nil(2)-module $\partial $. The category of nil(2)-modules can be considered to be a full subcategory of
the category of quadratic modules. It is a
reflective subcategory since there is a reflection functor given by Baues in
\cite{baus}. This is of course functorial and we can say that there is a
functor from the category of nil(2)-modules to that of quadratic modules. We
denote it by
\begin{equation*}
D:\mathsf{Nil}(2)\rightarrow \mathsf{Quad}.
\end{equation*}
Therefore, the left adjoint of $\Phi$ assigns to a nil(2)-module $\partial:M\rightarrow N$ the quadratic module written
\begin{equation*}
\xymatrix{ & C\otimes C\ar[dl]_{1} \ar[d]^{w} \\ L \ar[r]_{w} & M
\ar[r]_{\partial} & N. }
\end{equation*}
Our first main result is:
\begin{prop}\label{fibreqm}
The forgetful functor $\Phi:\mathsf{Quad}\rightarrow \mathsf{Nil}(2)$ is fibred and has a left adjoint.
\end{prop}
\begin{pf}
The left adjoint was given above. To prove that $\Phi$ is fibred, we give the  pullback construction of a quadratic module.

 Consider
the quadratic module
\begin{equation*}
\xymatrix{&C\otimes C\ar[d]_w\ar[dl]_{\omega}& \\\sigma:L \ar[r]_{\partial_2}&
M\ar[r]_{\partial_1}&N}
\end{equation*}
with a morphism of nil(2)-modules $u:=(u_1,u_0)$ given by the following commutative diagram:
$$
\xymatrix{R\ar[r]_{u_1}\ar[d]_{\partial_1'}&M\ar[d]^{\partial_1} \\ S\ar[r]_{u_0}&N.}
$$
We shall construct a quadratic module on the nil(2)-module $\partial_1':R\rightarrow S$.

Let
$$u^*(L)=\{(r,l):r\in
\ker \partial_1', u_1(r)=\partial_2(l)\}\subset R\times L.$$
We have a commutative
diagram
\begin{equation*}
\xymatrix{u^{*}(L)\ar[d]_{\partial_2'} \ar[r]^{\varphi}
&L\ar[d]^{\partial_2}\\
R\ar[r]_{u_1}\ar[d]_{\partial_1'}&M\ar[d]^{\partial_1} \\ S\ar[r]_{u_0}&N}
\end{equation*}
in which  $\varphi(r,l)=l$ and $\partial_2^{\prime }(r,l)=r$ for $(r,l)\in u^{*}(L).$

Let $C'=(R^{cr})^{ab}$. We can define the quadratic map $\omega'$ from $C'\otimes C'$ to $L$ as
$$\omega^{\prime }(\{r\}\otimes\{r^{\prime
}\})=(\langle r,r^{\prime }\rangle , \omega \{u_1(r)\}\otimes\{u_1(r')\})$$
where $\omega$ is the quadratic map of the quadratic module $\sigma$.   The action of $S$ on $u^*(L)$ is given by
$(r,l)^s=(r^s,l^{u_0(s)})$ for $s\in S$, $(r,l)\in u^*(L)$. Then since $
u_1(r^s)=u_1(r)^{u_0(s)}=\partial_2(l)^{u_0(s)}=\partial_2(l^{u_0(s)})$
and $\partial_1'$ is a pre-crossed module, for  $r\in \ker\partial_1'$ we have $r^s\in \ker\partial_1'$ and $(r^s,l^{u_0(s)})\in u^*(L)$.

We obtain that the diagram of homomorphisms of groups
\begin{equation*}
\xymatrix{&C'\otimes C'\ar[d]_{w'}\ar[dl]_{\omega'}& \\u^{*}(L)
\ar[r]_{\partial_2'}& R\ar[r]_{\partial_1'}&S}
\end{equation*}
is a quadratic module. (See [$\bullet^1]$ in Appendix)

It can be easily shown that $(\varphi,u_1,u_0)$ is a morphism between quadratic modules.
 Now we show that the morphism $(\varphi,u_1,u_0)$ is a cartesian morphism in the category of quadratic
 modules over the morphism $\Phi(\varphi,u_1,u_0)=(u_1,u_0)$.  Let $(\upsilon_1,\upsilon_0):\partial_1''\rightarrow
\partial_1'$ be morphism between nil(2)-modules given by the commutative diagram
$$
\xymatrix{K_1\ar[r]^{\upsilon_1}\ar[d]_{\partial_1''}&R\ar[d]^{\partial_1'}\\
K_0\ar[r]_{\upsilon_0}&S}
$$
and
\begin{equation*}
\xymatrix{&C''\otimes C''\ar[d]_{w''}\ar[dl]_{\omega''}& \\Z
\ar[r]_{\partial_2''}& K_1\ar[r]_{\partial_1''}&K_0}
\end{equation*}
be quadratic module with a morphism $(\theta, u_1^{\prime },u_0^{\prime })$ from the quadratic module $(\omega'',\partial_2'',\partial_1'')$ to the quadratic module $(\omega,\partial_2,\partial_1)$ together with $u_1\upsilon_1=u'_1$ and $u_0\upsilon_0=u'_0$.

We have the following commutative diagram
\begin{equation*}
\xymatrix{Z\ar[d]_{\partial''_2}\ar@{.>}[r]_\psi\ar@/^1.5pc/[rrr]^\theta
&u^*(L)\ar[d]_{\partial'_2}\ar[rr]_\varphi&&L\ar[d]^{\partial_2}\\
K_1\ar[d]_{\partial''_1}\ar@/^1.5pc/[rrr]_{u_1'}\ar[r]_{\upsilon_1}&R
\ar[d]_{\partial'_1}\ar[rr]_{u_1}&&M\ar[d]^{\partial_1}\\
K_0\ar@/^1.5pc/[rrr]_{u_0'}\ar[r]_{\upsilon_0}&S
\ar[rr]_{u_0}&&N.}
\end{equation*}
The necessary unique quadratic module morphism $(\psi,\upsilon_1,\upsilon_0)$ is given by
$\psi(z)=(\upsilon_1\partial''_2(z),\theta(z))$, for $z\in Z$. We show  $[\bullet^2]$ in Appendix that $(\varphi,u_1,u_0)$ is a morphism of quadratic modules.

Then the morphism $(\varphi,u_1,u_0)$ is a cartesian morphism in $\mathsf{Quad}$ over the morphism
$\Phi(\varphi,u_1,u_0)=u:=(u_1,u_0)$. The quadratic module
\begin{equation*}
\xymatrix{&C'\otimes C'\ar[d]_{w'}\ar[dl]_{\omega'}& \\u^{*}(L)
\ar[r]_{\partial_2'}& R\ar[r]_{\partial_1'}&S}
\end{equation*}
is called the \textit{pullback} of the quadratic module
$$
\xymatrix{&C\otimes C\ar[d]_w\ar[dl]_{\omega}& \\\sigma:L \ar[r]_{\partial_2}&
M\ar[r]_{\partial_1}&N}
$$
by the morphism $u:=(u_1,u_0)$ between nil(2)-modules.

\end{pf}

We may select a cartesian lifting of $\sigma$ along $u:=(u_1,u_0)$
$$
u^{\sigma}:u^*(\sigma)\rightarrow \sigma
$$
where $u^*(\sigma)$ is the pullback quadratic module
\begin{equation*}
\xymatrix{&C'\otimes C'\ar[d]_{w'}\ar[dl]_{\omega'}& \\u^{*}(L)
\ar[r]_{\partial_2'}& R\ar[r]_{\partial_1'}&S.}
\end{equation*}
If we fix the morphism $u:=(u_1,u_0)$  from $(\partial'_1:R\rightarrow S)$ to $(\partial_1:M\rightarrow N)$ in the category of nil(2)-modules, we obtain a re-indexing functor
$$
u^*:\mathsf{Quad}_{(M\rightarrow N)}\longrightarrow \mathsf{Quad}_{(R\rightarrow S)}
$$
defined on objects by $\sigma \mapsto u^*(\sigma)$ and where $\mathsf{Quad}_{(M\rightarrow N)}$ denotes the category of quadratic modules over the same  nil(2)-module $\partial_1:M\rightarrow N$.
In particular there is, for a fixed nil(2)-module $\partial:M\rightarrow N$, a subcategory $\mathsf{Quad}_{(M\rightarrow N)}$ of $\mathsf{Quad}$ which has
as objects those quadratic modules with $\partial:M\rightarrow N$
 as the `base'.

Using Proposition 2.5 of \cite{brown6}, for this re-indexing functor $u^*$, there is a bijection
$$
\mathsf{Quad}_{(R\rightarrow S)}(\kappa,u^*(\sigma))\cong \mathsf{Quad}_{(u_1,u_0)}(\kappa,\sigma)
$$
natural in $\kappa\in \mathsf{Quad}_{(R\rightarrow S)}$, $\sigma\in \mathsf{Quad}_{(M\rightarrow N)}$ and where $\mathsf{Quad}_{(u_1,u_0)}$ consists of those morphisms $(\alpha,u_1,u_0):\kappa\rightarrow \sigma$
 in $\mathsf{Quad}$ with  $\Phi(\alpha,u_1,u_0)=u:=(u_1,u_0)$.

\section{Quadratic Modules Cofibred Over Nil(2)-modules}
We give the construction of cofibration of $\Phi$.
Brown and Higgins in \cite{brown4} described the \textit{induced} crossed module $u_*(M)\rightarrow Q$ from a crossed module $M\rightarrow P$ by a morphism $u:P\rightarrow Q$
in the category of groups. That is, they proved that the forgetful functor from the category of crossed modules to the category of groups which
 sends $(M\rightarrow P)\mapsto P$ is a fibration and also a cofibration of categories.  We will extend  this result for the functor
$$
\Phi:\mathsf{Quad}\rightarrow \mathsf{Nil}(2).
$$
This gives that quadratic modules cofibred over nil(2)-modules and the notion of  \textit{induced quadratic module}. By a similar way, the notion of
induced 2-crossed module has been constructed in \cite{AAO}.
\begin{prop}\label{cofibreqm}
The forgetful functor $\Phi:\mathsf{Quad}\rightarrow \mathsf{Nil}(2)$ is cofibred.
\end{prop}
\begin{pf}
We prove this by a direct construction.
Consider the quadratic module
\begin{equation*}
\xymatrix{&C\otimes C\ar[d]_w\ar[dl]_{\omega}& \\L \ar[r]_{\partial_2}&
M\ar[r]_{\partial_1}&P}
\end{equation*}
with a morphism of nil(2)-modules $v:(v_1,v_0)$ from $(M\rightarrow P)$ to $(N\rightarrow Q)$ given by a commutative diagram
$$
\xymatrix{M\ar[r]_{v_1}\ar[d]_{\partial_1}&N\ar[d]^{\partial'_1} \\ P\ar[r]_{v_0}&Q.}
$$
We shall construct a quadratic module on the nil(2)-module $N\rightarrow Q$.

For this, consider the free product
$$
F(N\times L)\ast \langle C\otimes C\rangle
$$
where $\langle C\otimes C\rangle $ is the
free group generated by the elements of the form $\langle \{a\}\otimes
\{b\}\rangle $ for $a,b\in N$ and $F(N\times L)$ is the free group generated
by the set $N\times L$ and $C=(N^{cr})^{ab}$. Observe that action of $Q$ on $N$ induces an action on $\langle
C\otimes C\rangle $ by
$$(\{a\}\otimes \{b\})^{q}=(\{a{^{q}}
\}\otimes \{b{^{q}}\})$$
and induces an action on $F(N\times L)$ given by $(n,l)^{q}=(n^{q},l)$.
To get an induced quadratic module, factor $F(N\times L)\ast \langle C\otimes C\rangle$ out by the relations:
\begin{enumerate}
\item $( \{n^{-1}v_1 \partial _{2}ln\}\otimes \{b\})(\{b\}\otimes \{n^{-1}v_1 \partial _{2}ln\})=
(n,l)^{-1}(n,l)^{\partial _{1}^{\prime }b}$
\item $( \{n^{-1}v_1\partial _{2}ln\}\otimes \{n^{\prime
-1}v_1\partial _{2}l^{\prime }n^{\prime -1}\}=[(n,l),(n^{\prime
},l^{\prime })]$
\end{enumerate}
for $n,n',b\in N,l,l'\in L$.

We have a morphism $\partial _{2}^{\prime }:F(N\times L)\ast \langle C\otimes C\rangle \rightarrow N$ induced on $F(N\times L)$ by
$$(n,l)\mapsto n^{-1}(v_{1}\partial _{2}l)n$$
and given on $\langle C\otimes C\rangle$ by
$$\langle\{x\}\otimes \{y\}\rangle\mapsto x^{-1}y^{-1}xy^{\partial'_1x}$$
for $x,y\in N$. Then, the diagram
$$
\xymatrix{F(N\times L)*\langle C\otimes C\rangle\ar[r]^--{\partial'_{2}}& N\ar[r]^{\partial'_{1}}&Q}
$$
is a complex of homomorphisms of groups. (See [$\bullet^3$] Appendix)

 We define the quadratic map
$$
\omega':C\otimes C\longrightarrow F(N\times L)*\langle C\otimes C\rangle
$$
by $\omega ^{\prime }(\{x\}\otimes \{y\})=\langle \{x\}\otimes \{y\}\rangle$
for $x,y\in N.$

We define  $\psi :L\rightarrow F(N\times L)*\langle C\otimes C\rangle$ by $\psi(l)=(1,l)$  in the following diagram,
\begin{equation*}
\xymatrix{L\ar[d]_{\partial_2}\ar[r]^-{\psi}&F(N\times L)*\langle
C\otimes C\rangle\ar[d]^{\partial_2'}\\
M\ar[r]_{v_1}\ar[d]_{\partial_1}&N\ar[d]^{\partial_1'} \\
P\ar[r]_{v_0}&Q.}
\end{equation*}
We now wish to change the map $\psi :L\rightarrow F(N\times L)*\langle C\otimes C\rangle$ to make it operator morphism
in $\mathsf{Quad}$. For this we need to add in  $ F(N\times L)*\langle C\otimes C\rangle$ new relations:
\begin{enumerate}
\item $(1,\omega\{x\}\otimes \{y\})=\langle \{v_1(x)\}\otimes \{v_1(y)\} \rangle   $
\item $(1,l^p)=(v_0p,l) $
\end{enumerate}
for $p\in P, l\in L$ and $x,y\in M$. Then $(\psi,v_1,v_0)$ is a morphism of quadratic modules. (See [$\bullet^4$] Appendix)

Now, we show that the morphism $(\psi,v_1,v_0)$ in the category of quadratic modules is cocartesian over the morphism $\Phi(\psi,v_1,v_0)=v:=(v_1,v_0)$
 in the category of nil(2)-modules. For this, we use Proposition \ref{2} given by Brown and Sivera in \cite{brown6}.

Let
$$
\xymatrix{&C\otimes C\ar[d]_{w''}\ar[dl]_{\omega''}& \\X\ar[r]_{\partial''_2}&N \ar[r]_{\partial'_1}&Q}
$$
be a quadratic module in $\mathsf{Quad}_{(N\rightarrow Q)}$ with a morphism of quadratic modules $(\theta ^{\prime },v_{1},v_{0})$ given by the commutative
diagram:
\begin{equation*}
\xymatrix{L\ar[d]_{\partial_2} \ar[rr]^{\theta'}& &X\ar[d]^{\partial''_2}\\
M\ar[rr]_{v_1}\ar[d]_{\partial_1}&&N\ar[d]^{\partial_1'} \\
P\ar[rr]_{v_0}&&Q.}
\end{equation*}

Then, there is a unique morphism $\theta _{\ast }:F(N\times L)\ast \langle
C\otimes C\rangle \rightarrow X$ such that the diagram
\begin{equation*}
\xymatrix{L\ar[d]_{\partial_2}\ar[r]_-\psi\ar@/^1.5pc/[rr]^{\theta'} &F(N\times
L)*\langle C\otimes C\rangle\ar[d]_-{\partial'_2}\ar@{.>}[r]_-{\theta_*}&X\ar[d]^{\partial''_2}\\
M\ar[d]_{\partial_1}\ar[d]\ar[r]_{v_1}&N\ar[d]_{\partial'_1}\ar@{=}[r]_{id}&N\ar[d]^{\partial'_1}\\
P\ar[r]_{v_0}&Q\ar@{=}[r]_{id}&Q}
\end{equation*}
commutes.

The necessary unique morphism $\theta_*$ is defined as follows:
For generators $(n,l)$ of $F(N\times L)$, $\theta_{\ast }$ is given by
$\theta_{\ast }(n,l)=(\theta'(l))^n$. For generators $\langle \{x\}\otimes \{y\}\rangle $ of $\langle C\otimes C\rangle$,  $\theta_*$ is defined by
 $\theta_*\langle \{x\}\otimes \{y\}\rangle =\omega''(\{x\}\otimes \{y\})$
where $\omega''$ is the quadratic map from $C\otimes C$ to $X$. (See [$\bullet^5$] in Appendix)
 Using Proposition \ref{2}, the morphism $(\psi,v_1,v_0)$ becomes a cocartesian morphism in $\mathsf{Quad}$ over $\Phi(\psi,v_1,v_0)=v:=(v_1,v_0)$.
Clearly, the forgetful functor
$$
\Phi:\mathsf{Quad} \rightarrow \mathsf{Nil}(2)
$$
is cofibred.
\end{pf}

Thus, together with the relations in $F(N\times L)*\langle C\otimes C\rangle$ given above, the diagram
$$
\xymatrix{&C\otimes C\ar[d]_{w}\ar[dl]_{\omega'}& \\F(N\times L)*\langle C\otimes C\rangle \ar[r]_-{\partial'_2}&N \ar[r]_{\partial'_1}&Q}
$$
is an \textit{induced quadratic module} by the morphism $v:=(v_1,v_0)$ of nil(2)-modules.

Consequently, for a nil(2)-module morphism
$$v:=(v_1,v_0):(\xymatrix{M\ar[r]^{\partial_1}&P})\longrightarrow(\xymatrix{N\ar[r]^{\partial'_1}&Q})$$
we obtain a functor
$$
v_*:\mathsf{Quad}_{(M\rightarrow P)}\longrightarrow \mathsf{Quad}_{(N\rightarrow Q)}
$$
which sends a  quadratic module
$$
\xymatrix{L\ar[r]^{\partial_2}&M\ar[r]^{\partial_1}&P}
$$
to the induced quadratic module
$$
\xymatrix{&C\otimes C\ar[d]_{w}\ar[dl]_{\omega'}& \\v_*(L)\ar[r]_{(\partial_2)_{*}}&N \ar[r]_{\partial'_1}&Q}
$$
where $v_*(L)=F(N\times L)*\langle C\otimes C\rangle$ with the relations and $(\partial_2)_{*}=\partial'_2$.

Therefore we obtain the following result.
\begin{thm}
Let
$$v:=(v_1,v_0):(\xymatrix{M\ar[r]^{\partial_1}&P})\longrightarrow(\xymatrix{N\ar[r]^{\partial'_1}&Q})$$
be a morphism of nil(2)-modules. The functor
$$
v_*:\mathsf{Quad}_{(M\rightarrow P)}\longrightarrow \mathsf{Quad}_{(N\rightarrow Q)}
$$
is left adjoint to the pullback functor
$$
v^*:\mathsf{Quad}_{(N\rightarrow Q)}\longrightarrow\mathsf{Quad}_{(M\rightarrow P)}.
$$
\end{thm}

\begin{pf}
Let $$
\xymatrix{L\ar[r]^{\partial_2}&M\ar[r]^{\partial_1}&P}
$$
be a quadratic module. The diagram
$$
\xymatrix{L\ar@/^1.5pc/[drr]^{\partial_2}\ar[dr]^{\eta_L}\ar@/_1.5pc/[ddr]_{\psi}&&&\\
&v^{*}(v_*(L))\ar[r]_-{(\partial_2)_*^{*}}\ar[d]_{\varphi}&M\ar[r]^{\partial_1}\ar[d]^{v_1}&P\ar[d]^{v_0}\\
&v_*(L)\ar[r]_{(\partial_2)_*}&N\ar[r]_{\partial'_1}&Q}
$$
determines the morphism in $\mathsf{Quad}_{(M\rightarrow P)}$
$$
(\eta_L,id,id):(\xymatrix{L\ar[r]^{\partial_2}&M\ar[r]^{\partial_1}&P})\rightarrow (\xymatrix{v^*(v_*L)\ar[r]^{(\partial_2)_*^{*}}&M\ar[r]^{\partial_1}&P}).
$$
We show that this is a universal morphism.
Let
$$
\xymatrix{X\ar[r]^{\partial''_2}&N\ar[r]^{\partial'_1}&Q}
$$
be a quadratic module. Given a morphism
$$
(t,id,id):(\xymatrix{L\ar[r]^{\partial_2}&M\ar[r]^{\partial_1}&P})\rightarrow (\xymatrix{v^*(X)\ar[r]^{(\partial''_2)^{*}}&M\ar[r]^{\partial_1}&P})
$$
in $\mathsf{Quad}_{(M\rightarrow P)}$, we consider the composition morphism
$$
(\varphi t,v_1,v_0):(\xymatrix{L\ar[r]^{\partial_2}&M\ar[r]^{\partial_1}&P})\rightarrow (\xymatrix{X\ar[r]^{(\partial''_2)}&N\ar[r]^{\partial'_1}&Q}).
$$
By Proposition \ref{cofibreqm}, there exists only one morphism
$$
(\xymatrix{v_*(L)\ar[r]^{(\partial_2)_*}&N\ar[r]^{\partial'_1}&Q})\rightarrow (\xymatrix{X\ar[r]^{(\partial''_2)}&N\ar[r]^{\partial'_1}&Q})
$$
in $\mathsf{Quad}_{(N\rightarrow Q)}$, which commutes the following diagram:
$$
\xymatrix{\left(L\rightarrow M\rightarrow P\right)\ar[rr]^{(\psi,v_1,v_0)}\ar[dd]_{(\varphi t,v_1,v_0)}&&\left(v_*(L)\rightarrow N\rightarrow Q\right)\ar[ddll]^{((\varphi t)_*,id,id)}\\
\\
\left(X\rightarrow N\rightarrow Q\right)&&}
$$
Through the construction of $\eta_L$ and $v^*$ it is evident that:
$$
\xymatrix{(L,M,P)\ar@{=}[d]\ar[rr]^{(\eta_L,id,id)}&&v^{*}(v_*L)\ar[d]^-{v^*((\varphi t)_*,id,id)}\ar[r]^-{(\partial_2)_{*}^{*}}&(M\rightarrow P)\ar@{=}[d]\\
(L,M,P)\ar[rr]_-{(t,id,id)}&&v^{*}X\ar[r]_{(\partial''_2)^{*}}&(M\rightarrow P)}
$$
since:
$$
\xymatrix@!0{
L\ar@/_1.2pc/[dddd]_{\psi}\ar[ddr]_{\eta_L}\ar@/^2.7pc/[ddddddr]^{t}\ar@/^1.5pc/[ddrrrrr]^{\partial_2}\\
\\
  & v^*(v_*L)\ar[rrrr]_{(\partial_2)_{*}^{*}} \ar[dddd]\ar[ddl]_{\varphi}
      &  & & &(M\rightarrow P) \ar@{=}[dddd]\ar[ddl]^{v}        \\
      \\
  v^*L \ar[rrrr]\ar[dddd]_{((\varphi t)_*,id,id)}
      &  &  & &(N\rightarrow Q) \ar@{=}[dddd] \\
      \\
  & v^*(X) \ar[rrrr]^{(\partial''_2)^*}\ar[ddl]_{\varphi}
      &  & &   &(M\rightarrow P)\ar[ddl]^{v}               \\
      \\
  X \ar[rrrr]_{\partial''_2}
      &  &  &  & (N\rightarrow Q)    }
$$
Thus $v_*$ is left adjoint to $v^*$.

\end{pf}

\subsection{Another Presentation of Induced Quadratic Modules}

 Crossed squares were introduced in \cite{walery}. Quadratic modules are related to crossed squares,
 and induced crossed squares have been studied in \cite{brown9, ellis} and also in   Appendix B4 of \cite{brown5}.
The following method used in \cite{brown6} gives another view of a presentation of the
induced crossed square, and which is applied to free crossed squares in \cite{ellis}.

In Theorem \ref{3}, by taking $\mathsf{X}=\mathsf{Quad}$ and $\mathsf{B}=\mathsf{Nil}(2)$,
we can give another view of a description of the induced quadratic module.

We know that the functor $\Phi:\mathsf{Quad} \rightarrow \mathsf{Nil}(2)$ is a fibration and has a left adjoint
$$
D:\mathsf{Nil}(2)\rightarrow \mathsf{Quad}.
$$
Recall from the previous section that this left adjoint $D$ assigns  to a nil(2)-module $\partial:M\rightarrow P$ the quadratic module written
$$
\xymatrix{ & C\otimes C\ar[dl]_{1} \ar[d]^{w} \\ L \ar[r]_{w} & M
\ar[r]_{\partial} & P }
$$
where $L=C\otimes C$. We denote this associated quadratic module by
$$D(M\rightarrow P)=\xymatrix{C\otimes C\ar[r]^{w}&M\ar[r]&P}.$$

 We know from \cite{baus} that the category of quadratic modules admits
pushouts. Then $\Phi$ is a fibration of categories and  also a cofibration.
Thus we have a notion of \textit{induced quadratic module}, as follows:

Given a quadratic module
\begin{equation*}
\xymatrix{& C\otimes C\ar[d]^{w}\ar[dl]_{\omega} &\\ L \ar[r]_{\delta} &
M\ar[r]_{\partial}&P }
\end{equation*}
and a morphism of nil(2)-modules
$$
\xymatrix{M\ar[r]^{v_1}\ar[d]_{\partial}&N\ar[d]^{\partial'} \\ P\ar[r]_{v_0}&Q}
$$
we get an \textit{induced quadratic module}
$$
\xymatrix{& C'\otimes C'\ar[d]^{w}\ar[dl]_{\omega'} &\\ T \ar[r]_{\delta'} &
N\ar[r]_{\partial'}&Q }
$$
which according to Theorem \ref{3} is given by a pushout in the category of quadratic modules of the form
\begin{equation*}
\xymatrix{(C\otimes C\rightarrow M\rightarrow P)\ar[dd]_{(u,id,id)}\ar[rr]^{(v_1^{cr^{ab}},v_1,v_0)}&&(C'\otimes C'\rightarrow N\rightarrow Q)
\ar[dd]^{(v,id,id)}\\
\\
(L\rightarrow M\rightarrow P)\ar[rr]_{(\psi,\upsilon_1,\upsilon_0)}&&(T\rightarrow N\rightarrow Q).}
\end{equation*}
This gives another view of a presentation of the induced quadratic module given above.

\section{Applications to Reduced Quadratic Modules}

Recall that a reduced quadratic module is a quadratic module
\begin{equation*}
\xymatrix{& C\otimes C\ar[d]^{w}\ar[dl]_{\omega} &\\ L \ar[r]_{\delta} &
M\ar[r]_{\partial}&N }
\end{equation*}
in which the group N is a trivial group.
More clearly, a reduced quadratic module $(\omega,\delta)$ is a diagram
\begin{equation*}
\xymatrix{ M^{ab}\otimes M^{ab}\ar[r]^-{\omega} & L \ar[r]^{\delta} & M }
\end{equation*}
of homomorphisms between groups such that the following conditions hold:%
\newline
$(i)$ The group $M$ is a nil(2)-group and the quotient map $%
M\twoheadrightarrow M^{ab}$ to the \emph{abelianization} $M^{ab}$ of $M$ is
denoted by $x\mapsto \overline{x}$.\newline
$(ii)$ For $x,y\in M$, $\delta\omega(\overline x\otimes\overline y)=[x,y].$\newline
$(iii)$ For $a\in L, x\in M$, $\omega(\overline{\delta a}\otimes\overline x)(\overline x\otimes\overline{\delta a})=1.$\newline
$(iv)$ For $a,b\in L$, $\omega(\overline{\delta a}\otimes \overline{\delta b})=[a,b].$

We denote the category of reduced quadratic modules by $\mathsf{RQM}$. Since in a reduced quadratic module
\begin{equation*}
\xymatrix{ C\otimes C\ar[r]^-{\omega} & L \ar[r]^{\delta} & M }
\end{equation*}
the group $M$ is a nil(2)-group, we have a forgetful functor
$$\Phi_R:\mathsf{RQM}\longrightarrow \mathsf{NilGrp}(2)$$
from the category of reduced quadratic modules to that of nil(2)-groups which sends
\begin{equation*}
\xymatrix{ (C\otimes C\ar[r]^-{\omega} & L \ar[r]^{\delta} & M)\ar[r]&M. }
\end{equation*}
Using the construction method given in the previous section, we can give the following  proposition.
\begin{prop}
The functor $\Phi_R:\mathsf{RQM}\rightarrow \mathsf{NilGrp}(2)$ is a bifibration.
\end{prop}

We now compare the universal properties defining the induced reduced quadratic module and the free reduced quadratic module on a map.
Using the induced reduced quadratic module constructed above, we get an alternative description of the free reduced quadratic module.

\textbf{Free Reduced Quadratic Modules}

In \cite{au},  using the suspension functor given in \cite{Muro}, the construction of a free reduced quadratic module was given.
 We recall the definition of free reduced quadratic modules. Let
\begin{equation*}
\xymatrix{ C\otimes C\ar[r]^-{\omega} & L \ar[r]^{\delta} & M }
\end{equation*}
be a reduced quadratic module, $Y$ be a set and $\upsilon: Y\rightarrow L$ be a function. Then this reduced quadratic module is said to be a free reduced quadratic module on function $\upsilon: Y\rightarrow L$ if for any reduced quadratic module
\begin{equation*}
\xymatrix{ C'\otimes C'\ar[r]^-{\overline\omega} & L' \ar[r]^{\overline\delta} & M }
\end{equation*}
and a function $\upsilon': Y\rightarrow L'$ such that $\delta\upsilon=\overline\delta\upsilon'$, there is a unique morphism $\Phi: L\rightarrow L'$ such that $\overline\delta\Phi=\delta.$
This situation can be illustrated as follows:
$$\xymatrix{C\otimes C\ar[rr]^{\omega}\ar[dd]_{\Phi^{*}}&&L \ar[dd]_{\Phi}\ar[rr]^{\delta}&&M\ar@{=}[dd]\\
&&&Y \ar[ul]_{\upsilon}\ar[dl]^{\upsilon'}&\\
C'\otimes C'\ar[rr]_{\omega'}&&L'\ar[rr]_{\overline{\delta}}&&M }
$$

\begin{prop}
Let P be a nil(2)-group and $\{\theta_r:r\in R\}$ be an indexed family of elements of P, or equivalently, a map $\theta:R\rightarrow P$.
Let F be the free group generated by R and $F^{nil}$ be the nilization of F with the quotient map $g:F\rightarrow F^{nil}, r\mapsto \overline r$.
Let $f^{nil}:F^{nil}\rightarrow P$ be the homomorphism of groups such that $f^{nil}(\overline r)=\theta_r \in P$. We have the following reduced quadratic module
\begin{equation*}
\xymatrix{ (\omega,1)=C\otimes C\ar[r]^-{\omega} & F^{nil} \ar[r]^{1_{F^{nil}}} & F^{nil}}
\end{equation*}
where $\omega$ is given by $\omega\{\overline r\}\otimes\{\overline {r'}\}=[\overline r,\overline {r'}]$ and where $\overline{r} \in F^{nil} $,
 $\{\overline {r}\} \in C$. Then the reduced quadratic module
\begin{equation*}
\xymatrix{ C'\otimes C'\ar[r] & {F(P\times F^{nil})*\langle C'\otimes C'\rangle={f_*}^{nil}(P)}\ar[r]^-{\overline\delta} & P }
\end{equation*}
induced from $(\omega,1)$ by $f^{nil}$ is the free reduced quadratic module on $\{(1,\overline r):r\in R\}$.
\end{prop}
\begin{pf}
First, consider the following diagram
\begin{equation*}
\xymatrix{{F^{nil}}\ar[dd]_{1_{F^{nil}}}\ar[rr]^{\phi}\ar@/^1.3pc/[dddrrr]^---------------{h^{nil}}&&{f_*}^{nil}(P)\ar[dd]_{\overline\delta}\ar@/^1.3pc/[dddr]^{\phi_*}&\\
&&&&\\
F^{nil}\ar[rr]_{f^{nil}}&&P&\\
&R\ar[ul]^{inc}\ar[ur]_{\theta}\ar[rr]_{\theta'}&&N\ar[ul]}
\end{equation*}
in which $\phi(\overline r)=(1,\overline r)$.
We shall check that the data in the free and induced constructions are equivalent. The data in the induced construction are a reduced quadratic module
\begin{equation*}
\xymatrix{ (\omega'',\partial):C'\otimes C'\ar[r]^-{\omega''} & N \ar[r]^{\partial} & P }
\end{equation*}
and a morphism of $(h^{nil},f^{nil}):(\omega',1)\rightarrow(\omega'',\partial)$. The data in the free reduced quadratic module
construction are a reduced quadratic module
\begin{equation*}
\xymatrix{ (\omega'',\partial):C'\otimes C'\ar[r]^-{\omega''} & N \ar[r]^{\partial} & P }
\end{equation*}
and a map $\theta':R\rightarrow N$ with $\partial\theta'=\theta$. Since $F$ is the free group on $R$, the map
$\theta'$ is equivalent to a homomorphism $h:F\rightarrow N$ lifting $\theta$ (i.e. $h(r)=\theta'(r)$).
We have a morphism from $h$, $h^{nil}:F^{nil}\rightarrow N$  by $h^{nil}(\overline r)=h(r)=\theta'(r))$. Moreover $h^{nil}$ satisfies
\begin{eqnarray*}
h^{nil}\omega\{\overline r\}\otimes\{\overline {r'}\} &=&h^{nil} [\overline r,\overline {r'}] \\
 &=& h^{nil}([\overline r,\overline {r'}]) \\
 &=& [\theta'(r),\theta' (r')] \\
 &=&  \omega''\{\partial\theta'r\}\otimes\{\partial\theta'{r'}\}\\
&=&  \omega''\{\theta_r\}\otimes\{\theta_{r'}\}\\
&=&  \omega''\{f^{nil}(\overline r)\}\otimes\{f^{nil}(\overline {r'})\}
\end{eqnarray*}
for all $r,r'\in R$. So $(h^{nil},f^{nil})$ is a morphism of reduced quadratic modules. Thus the data in both cases are equivalent.
\end{pf}

Now, we develop an example of a finite reduced quadratic module to see what it looks like and how the induced construction behave on it.
First, we give the following result.
\begin{prop}
Let $\phi:M\rightarrow N$ be a monomorphism of groups and let
$$\xymatrix{ C\otimes C\ar[r]^-{\omega} & L \ar[r]^{\partial} & M }$$
 be a reduced quadratic module and $T$ be a right transversal of $\phi(M)$ in $N$ and $(L\overrightarrow*T)*\langle C'\otimes C'\rangle$ be the free product
 of $(L\overrightarrow*T)$ and $\langle C'\otimes C'\rangle$ where
$(L\overrightarrow*T)$ is the free product of groups $L_t, t\in T$, with elements $(l,t),l\in L$  isomorphic to $L$ under the map
$(l,t)\rightarrow l$ and $\langle C'\otimes C'\rangle$ is the free group generated by the set $ C'\otimes C'$.
If $n\in N$ acts on $(L\overrightarrow*T)$ by the rule ${}^n(l,t)={}^u(l,m)$ where $m\in M, u\in T,$ $nt=u\phi(m)$, then
\begin{equation*}
\xymatrix{ C'\otimes C'\ar[r]^-{\omega} & \phi_*(L)=(L\overrightarrow*T)*\langle C'\otimes C'\rangle/S \ar[r]^-{\beta} & Q }
\end{equation*}
is a reduced quadratic module with the map $\omega(\overline a\otimes \overline b)=\langle\overline a\otimes \overline b\rangle S$ for
$\overline a, \overline b\in N^{ab}$, where $S$ is the normal closure in $\phi_*(L)$ of elements
$$\langle\overline{\beta(l,t)}\otimes \overline{\beta(l',t')} \rangle=[(l,t),(l',t')]$$
$$\langle\overline{\beta(l,t)}\otimes \overline{n}\rangle\langle \overline{n}\otimes \overline{\beta(l,t)}\rangle=1$$
$$\omega\{m\}\otimes\{m'\}=\langle\overline{\phi(m)}\otimes \overline {\phi(m')}\rangle$$
for $l \in L, t \in T, m,m' \in M, n\in N$.
\end{prop}

\begin{ex}(The Dihedral Reduced Quadratic Module)
Let
 \begin{equation*}
\xymatrix{ P^{ab}\otimes P^{ab}\ar[r]^-{\omega} & M \ar[r]^{\partial} & P }
\end{equation*}
be a reduced quadratic module and $i:P\rightarrow Q$ be a homomorphism. Suppose that $Q=D_4$ with the
 presentation $\langle a,b:a^4=b^2=abab=1\rangle$ is the dihedral group of order 8. Let $M=P$
  be the cyclic subgroup of $D_4$ of order 2 generated by $b$. Since $Z(D_4)=\{1,a^2\}$, we have $D_4/Z(D_4)\cong Klein 4-group=C.$
 Since $C$ is an Abelian group, the dihedral group $D_4$ is a nil(2)-group. Let $C_4=\{0,1,2,3\}$ be the cyclic group of order 4.
  A right transversal $T$ of $C_2=\{1,b\}$ in $D_4$ is given by the elements $1,x,x^2,x^3$,
  since $D_4/C_2=\{C_2,aC_2,a^2C_2,a^3C_2\}$. We have $i_*(C_2)=(C_2\overrightarrow*T)*\langle C\otimes C\rangle$ where
 $i:C_2\rightarrow D_4$ is
 the inclusion map and $C=(D_4)^{ab}$ is Klein 4-group.
  Hence the generators of $(C_2\overrightarrow*T)$ are in the forms $r_i=(b,a^i)$ with
  the relations $(r_i)^2=1, i=0,1,2,3$ and $\langle C\otimes C\rangle$ is the free group generated by
  elements of the forms $\langle ba^i\otimes ba^j\rangle$, $i,j=0,1,2$ with the relations
$$1. \langle ba^{2i}\otimes ba^{2j}\rangle= a^{2i+2j}\in \langle a^2\rangle$$
$$2. \langle ba^{2i}\otimes ba^j\rangle\langle ba^j\otimes ba^{2i}\rangle=1.$$
We have a morphism
$$\delta:(C_2\overrightarrow*T)*\langle C\otimes C\rangle\longrightarrow D_4$$
defined on generators $(b,a^i)$ by $\delta(b,a^i)=ba^{2i}$ and on generators $\langle ba^i\otimes ba^j\rangle$  by
$$\delta(\langle ba^i\otimes ba^j\rangle)=a^{2i+2j}.$$
Thus we have
 \begin{equation*}
 \begin{array}{rlcl}
   \delta(b\otimes a^2)=a^4=1= \delta( a^2\otimes b),& \delta(b\otimes a)=a^2=[b,a]=[a,b]=\delta(a\otimes b),\\
  \delta(b\otimes b)=1=[b,b], & \delta(a^2 b\otimes a)=1=[a^2,a],\\
  \delta(ba\otimes b)=[ba,b]=a^2, &\delta(ba\otimes ba)=a^4=1=[ba,ba],  \\
   \delta(b\otimes ba)=a^2=[b,ba].&
\end{array}
 \end{equation*}
Thus we have a reduced quadratic module
\begin{equation*}
\xymatrix{ C\otimes C\ar[r]^-{\omega} &(C_2\overrightarrow*T)*\langle C\otimes C\rangle \ar[r]^-{\delta} & D_4}
\end{equation*}
where $\omega\{ba^i\}\otimes\{ba^j\}=\langle ba^i\otimes ba^j\rangle$ for $i,j=0,1,2$. Moreover we define $u=r_0r_1=a,\upsilon=r_0=b$ and
$D_4'=\langle u,\upsilon: u^4=\upsilon^2=u\upsilon u \upsilon=1\rangle$ the another copy of $D_4$. Then we obtain $\delta(u)=a^2$ and $\delta(\upsilon)=b$ and
\begin{equation*}
\xymatrix{ C\otimes C\ar[r]^-{\omega} & D_4'*\langle C\otimes C\rangle \ar[r]^-{\delta} & D_4 }
\end{equation*}
is a dihedral reduced quadratic module, where $C$ is the Klein 4-group.
\end{ex}

\underline{REMARK:}  Ellis in his paper \cite{ellis} has stated that crossed squares have a geometric interpretation in terms of relative homotopy groups, but no such interpretation is 
available for quadratic modules. He also defined a functor $\rho$ from the category of 3-dimensional reduced CW-spaces to the category of crossed squares and showed that
$\rho(X)$ is a free crossed square. Furthermore in page 106 of his work  \cite{ellis}, obtained a free quadratic module from a free crossed square. Thus
 Brown and Sivera results for induced crossed squares have immediate topological applications.
Of course this description of induced structures for quadratic modules was very useful in situations where a van
 Kampen theorem gave pushouts arising from topology. 
 The problem with the quadratic module area is that there seems to be no van Kampen type theorem, 
 because there is no direct homotopical functor.

\section{Appendix}

\textbf{The Proof of Proposition \ref{fibreqm}}

$[\bullet^1$] We show that the following diagram
\begin{equation*}
\xymatrix{&C'\otimes C'\ar[d]_{w'}\ar[dl]_{\omega'}& \\u^{*}(L)
\ar[r]_{\partial_2'}& R\ar[r]_{\partial_1'}&S}
\end{equation*}
is a quadratic module.

$\mathsf{QM1}$. It is already known that $\partial':R\rightarrow S$ is a nil(2)-module and, since $r\in \ker \partial'_1$, we obtain
$\partial_1^{\prime }\partial_2^{\prime }(r,l)=\partial_1^{\prime }(r)=1.$ That is,
$$
\xymatrix{u^{*}(L)\ar[r]^{\partial_2'}&R\ar[r]^{\partial'_1}&S}
$$
is a complex of homomorphisms of groups.

$\mathsf{QM2}$. For all $r,r^{\prime }\in R$, we obtain
\begin{equation*}
\partial_2^{\prime }\omega^{\prime }(\{r\}\otimes \{r^{\prime
}\})=\partial_2^{\prime }(\langle r,r^{\prime }\rangle, \omega\{u_1r\}\otimes
\{u_1r'\})=\langle r,r'\rangle=w'(\{r\}\otimes \{r'\}).
\end{equation*}

$\mathsf{QM3}$. For $(r,l)\in u^*(L)$, $r^{\prime }\in R$, we obtain

\begin{multline*}
\omega^{\prime }(\{\partial_2^{\prime }(r,l)\}\otimes \{(r^{\prime
})\}\{r^{\prime }\}\otimes \{\partial_2^{\prime }(r,l)\})\\
\begin{aligned}
&=\omega^{\prime }(\{r\}\otimes \{r'\}\{r'\}\otimes \{r\})\\
&=( \langle r,r'\rangle,\omega\{u_1r\}\otimes \{u_1r' \})( \langle r',r\rangle,\omega\{u_1r'\}\otimes \{u_1r\}) \\
&=\left((r^{-1}r^{\prime -1}rr^{\prime \partial_1^{\prime }r})((r')^{-1}r^{-1}(r') r^{\partial_1'r'}),\omega\{\partial_2 l\}\otimes \{u_1r' \}\{u_1r'\}\otimes \{\partial_2 l\}\right)\\
&=\left((r^{-1}r^{\prime -1}rr^{\prime \partial_1^{\prime }r})((r')^{-1}r^{-1}(r') r^{\partial_1'r'}), l^{-1}l^{\partial_1 u_1(r')}\right) (\text{ since }u_1(r)=\partial_2l) \\
&=(r^{-1}r^{\partial_1^{\prime }r^{\prime }},l^{-1}l^{\partial_1 u_1(r')})\qquad (\text{ since } r\in \ker\partial_1') \\
&=(r,l)^{-1}(r^{\partial_1'(r')},l^{u_0\partial_1'r'})\\
&=(r,l)^{-1}(r,l)^{\partial_1'r'}
\end{aligned}
\end{multline*}

$\mathsf{QM4.}$ For $(r,l), (r^{\prime },l^{\prime })\in u^*(L)$, we get
\begin{align*}
\omega^{\prime }(\{\partial_2^{\prime }(r,l)\}\otimes \{\partial_2^{\prime
}(r,l)\})&=( \langle r,r^{\prime }\rangle,\omega\{u_1r\}\otimes \{u_1r^{\prime }\})
\\
&=(r^{-1}r^{\prime -1}rr^{\prime \partial_1^{\prime
}r},\omega\{\partial_2l\}\otimes \{\partial_2l^{\prime }\}) \\
&=(r^{-1}r^{\prime -1}rr^{\prime },[l,l^{\prime }])\quad (\text{ since } r\in \ker\partial_1') \\
&=([r,r^{\prime }],[l,l^{\prime }]) \\
&=[(r,l),(r^{\prime },l^{\prime })].
\end{align*}

$[\bullet^2]$ We show that $(\varphi,u_1,u_0)$ is a morphism  between quadratic modules.
Since
$$u_1(\upsilon_1\partial''_2(z))=u'_1(\partial''_2 z)=\partial_2(\theta(z))$$
and $\partial_1^{\prime
}(\upsilon_1\partial_2^{\prime \prime }z)=\upsilon_0\partial_1^{\prime
\prime }\partial_2^{\prime \prime }z=1$, we have $\upsilon_1\partial''_2(z)\in\ker\partial'_1$, and  $(\upsilon_1\partial''_2(z),\theta(z))\in u^*(L)$.
We obtain
$$
\partial'_2\psi(z)=\partial'_2(\upsilon_1\partial''_2(z),\theta z)=\upsilon_1\partial''_2(z)
$$
for $z\in Z$ and
\begin{align*}
\psi(\omega''\{x\}\otimes \{y\})&=(\upsilon_1(\partial''_2\omega''\{x\}\otimes \{y\}),\theta(\omega''\{x\}\otimes \{y\}))\\
&=(\upsilon_1(\langle x,y\rangle),\omega\{u'_1(x)\}\otimes \{u'_1(y)\})\\
&=(\langle \upsilon_1x,\upsilon_1y\rangle, \omega\{u_1(\upsilon_1x)\}\otimes \{u_1(\upsilon_1y)\})\\
&=\omega'(\{\upsilon_1x\}\otimes \{\upsilon_1y\})
\end{align*}
for $\{x\}\otimes \{y\}\in C''\otimes C'' $. $\Box$

\textbf{The Proof of Proposition \ref{cofibreqm}}

$[\bullet^3]$ We have for $(n,l)\in F(N\times L)$,
\begin{align*}
\partial'_1(\partial'_2)(n,l)&=\partial'_1(n^{-1}v_1\partial_2(l)n)\\
&=\partial'_1(n^{-1})\partial'_1 v_1\partial_2(l)\partial'_1n\\
&=\partial'_1(n^{-1}) v_0(\partial_1\partial_2l)\partial'_1n  \quad (\text{ since } \partial'_1v_1=v_0\partial_1)\\
&=\partial'_1(n^{-1})v_0(1)\partial'_1n  \quad (\text{ since } \partial_1\partial_2l=1)\\
&=1
\end{align*}
and  for $\langle\{x\}\otimes \{y\}\rangle\in \langle C\otimes C\rangle$,
$$
\partial'_1(\partial'_2 \langle\{x\}\otimes \{y\}\rangle)=\partial'_1(x^{-1}y^{-1}xy^{\partial'_1x}) =1.
$$
 Thus the diagram
$$
\xymatrix{F(N\times L)*\langle C\otimes C\rangle\ar[r]^--{\partial'_{2}}& N\ar[r]^{\partial'_{1}}&Q}
$$
is a complex of homomorphisms of groups.

$[\bullet^4]$ We show that  $(\psi,v_1,v_0)$ is a morphism of quadratic modules. We obtain $\partial'_2\psi(l)=\partial'_2(1,l)=v_1\partial_2(l)$ for $l\in L$ and
$$
\psi(\omega\{x\}\otimes \{y\})=(1,\omega\{x\}\otimes \{y\})=\langle \{v_1(x)\}\otimes \{v_1(y)\} \rangle =\omega'(\{v_1(x)\}\otimes \{v_1(y)\})
$$
for $x,y\in M$ and
$$\psi(l^p)=(1,l^p)=(v_0p,l)=(1,l)^{v_0p}=\psi(l)^{v_0p}$$
for $p\in P$.

$[\bullet^5]$ We show that $(\theta_*,id,id)$ is a morphism of quadratic modules. We obtain for $n\in N$ and $l\in L$
$$
\theta_{\ast}((n,l)^{q})=\theta _{\ast }(n^q,l)=(\theta'(l))^{n^{q}}=(\theta'(l)^{n})^{id(q)}=(\theta_*(n,l))^{id(q)}
$$
 and for $x,y\in N$,
$$
\theta_{\ast }(\langle \{x\}\otimes \{y\}\rangle^{q^{\prime }})=\theta _{\ast }\langle \{x^{q^{\prime
}}\}\otimes \{y^{q^{\prime }}\}\rangle =\omega ^{\prime \prime
}(\{x^{q^{\prime }}\}\otimes \{y^{q^{\prime }}\})=\omega ^{\prime
\prime }(\{x\}\otimes \{y\})^{q^{\prime }}=(\theta _{\ast }\langle
\{x\}\otimes \{y\}\rangle )^{id(q^{\prime })}$$
Further, we obtain
\begin{align*}
\partial''_2(\theta_*(n,l))&=\partial''_2(\theta'(l)^n)\\
&=n^{-1}(\partial''_2(\theta'(l)))n\\
&=n^{-1}(v_1\partial_2(l))n\\
&=id(\partial'_2(n,l))
\intertext{and}
\partial''_2(\theta_*\langle\{x\}\otimes \{y\}\rangle)&=\partial''_2(\omega''\{x\}\otimes \{y\})\\
&=w(\{x\}\otimes \{y\})\\
&=x^{-1}y^{-1}xy^{\partial'_1x}\\
&=\partial'_2\langle\{x\}\otimes \{y\}\rangle
\end{align*}
for $(n,l)\in F(N\times L)$ and $\langle\{x\}\otimes \{y\}\rangle\in \langle C\otimes C\rangle$. Therefore, $(\theta _{\ast },id,id)$
is the unique morphism in $\mathsf{Quad}_{(N\rightarrow Q)}$. $\Box$

\vspace{1cm}
$
\begin{array}{lll}
\text{Hasan AT\.{I}K} & \qquad \text{Erdal ULUALAN} &  \\
\text{\.{I}stanbul Medeniyet University} & \qquad \text{Dumlup\i nar
University} &  \\
\text{Science Faculty} & \qquad \text{Science and Art Faculty} &  \\
\text{Mathematics Department} & \qquad\text{Mathematics Department} &  \\
\text{\.{I}stanbul, TURKEY} & \qquad \text{K\"{u}tahya, TURKEY} &  \\
\text{\url{hasanatik@yahoo.com}} & \qquad \text{\url{eulualan@gmail.com}, \url{erdal.ulualan@dpu.edu.tr}} &  \\
&  &
\end{array}
$

\end{document}